\documentclass[12pt]{amsart}
\usepackage{amsmath}
\usepackage{amsthm}
\usepackage{enumitem}
\usepackage[left=1.00in, right=1.00in, top=1.00in, bottom=1.00in]{geometry}
\usepackage{cleveref}
\numberwithin{equation}{section}
\newtheorem{thm}{Theorem}[section]
\newtheorem{cor}[thm]{Corollary}
\newtheorem{ex}[thm]{Example}

\renewcommand{\Re}{\operatorname{Re}}
\renewcommand{\Im}{\operatorname{Im}}
\title{An Algorithm for Numerically Inverting the Modular $j$-function}
\author{Ethan Alwaise}
\address{Department of Mathematics, UCLA,
Los Angeles, CA 90095}
\email{ealwais@ucla.edu}
\begin{document}
\maketitle
\begin{abstract}
The modular $j$-function is a bijective map from $X_0(1) \setminus \{\infty\}$ to $\mathbb{C}$. A natural question is to describe the inverse map. Gauss offered a solution to the inverse problem in terms of the arithmetic-geometric mean. This method relies on an elliptic curve model and the Gaussian hypergeometric series. Here we use the theory of polar harmonic Maass forms to solve the inverse problem by directly examining the Fourier expansion of the weight $2$ polar harmonic Maass form obtained by specializing the logarithmic derivative of the denominator formula for the Monster Lie algebra.
\end{abstract}

\section{Introduction and Statement of Results}

Let $j(\tau)$ be the $\textrm{SL}_2(\mathbb{Z})$-modular function defined by
$$j(\tau) := \frac{E_4(\tau)^3}{\Delta(\tau)} = \sum_{n=-1}^{\infty}c(n)e^{2\pi in\tau} = e^{-2\pi i\tau} + 744 + 196884e^{2\pi i\tau} + \cdots,$$
where $E_{2k}(\tau)$ is the weight $2k$ Eisenstein series and $\Delta(\tau) := (E_4(\tau)^3 - E_6(\tau)^2)/1728$ is the modular discriminant. It is well known that $j(\tau)$ parametrizes isomorphism classes of elliptic curves over $\mathbb{C}$ and gives a bijective map from the fundamental domain $X_0(1) \setminus \{\infty\}$ to $\mathbb{C}$. A natural question is to ask for a description of the inverse map.

Gauss offered a solution to the inverse problem in terms of the arithmetic-geometric mean (AGM) by making use of the theory of elliptic functions. The elliptic Weierstrass $\wp$-function satisfies the differential equation
$$(\wp'(z))^2 = (\wp(z))^3 - g_2\wp(z) - g_3,$$
where $g_2 := 60E_4(\tau)$ and $g_3 := 140E_6(\tau)$ are the \textit{elliptic invariants}. The cubic equation above defines an elliptic curve over $\mathbb{C}$ whose $j$-invariant is equal to $j(\tau)$. Given $\alpha \in \mathbb{C}$, one can find $\tau \in X_0(1) \setminus \{\infty\}$ such that $j(\tau) = \alpha$ by producing a elliptic curve model of the above form with $j$-invariant $\alpha$. Then $\tau$ is given by the ratio of the fundamental periods $\omega_1$ and $\omega_2$ of the associated $\wp$-function. The theory of elliptic functions tells us that the inverse of the $\wp$-function is an elliptic integral. Using this fact, one can show that $\omega_1$ and $\omega_2$ are given in terms of so-called period integrals. Gauss showed that the period integrals are left unchanged by replacing certain parameters with their arithmetic and geometric means. By passing to the limit, he was able to evaluate the period integrals in terms of the AGM. The AGM can then be numerically evaluated using the Gaussian hypergeometric series.

It is natural to ask for a solution to the inverse problem that relies only on the properties of $j(\tau)$ as a function and makes no reference to an elliptic curve model. The theory of polar harmonic Maass forms offers such a solution. As it turns out, the logarithmic derivatives of meromorphic modular forms are polar harmonic Maass forms, as was shown by Bringmann et al. in \cite{divmf}. The inverse problem can then be reformulated in terms of locating the pole of the logarithmic derivative of $j(\tau) - \alpha$. This can be done using the asymptotic formula for the Fourier coefficients of such polar harmonic Maass forms offered in \cite{divmf}. Using their work, we prove the following theorem, which can be found in the M.S. thesis \cite{MSBS} of the author:

\begin{thm}\label{mainthm}
Let $\alpha \in \mathbb{C}$ and let $z \in X_0(1) \setminus \{\infty\}$ such that $j(z) = \alpha$. Define
$$H_z(\tau) := -\frac{1}{2\pi i}\frac{j'(\tau)}{j(\tau) - \alpha} = \sum_{n=0}^{\infty}a(n)e^{2\pi in\tau}.$$
Write $z = x + iy$. Then $y$ is given by
$$y = \lim_{n \to \infty}\frac{\log\vert a(n) \vert}{2\pi n}.$$
If $\alpha = 0$, then $x = -\tfrac{1}{2}$. If $\alpha \neq 0$, let
$$c(n) = \begin{cases}
\Re(a(n))e^{-2\pi ny_0} &\text{if } \lim_{n\to\infty}\vert a(n) \vert e^{-2\pi ny} = 1, \\
\frac{1}{2}a(n)e^{-2\pi ny_0} &\text{otherwise},
\end{cases}$$
where $y_0 \approx y$ is obtained from the $b(n)$. Let $w_n = \cos^{-1}(c(n))$. Then an approximation for $x$ is given by one of the following formulas:
\begin{equation*}
\begin{aligned}
x &\approx \pm\frac{1}{2\pi}(w_n \pm w_{n-1}) \\
x &\approx \pm\frac{1}{2\pi}(w_n + w_{n-1} - 2\pi).
\end{aligned}
\end{equation*}

\noindent Two remarks:

\begin{enumerate}
\item \emph{There is some ambiguity in the value of $x$ in the above theorem. However, it is not difficult to determine the correct value of $x$ by resubstituting the possible values into $j(\tau)$.}
\item \emph{It would be interesting to study the convergence of the above algorithm in detail.}
\end{enumerate}
\end{thm}

This paper is organized as follows: In Section 2, we briefly recall the basic facts about harmonic Maass forms. In Section 3, we use the work of Bringmann et al. in \cite{divmf} to prove \Cref{mainthm}. We conclude with Section 4, where we offer some examples of \Cref{mainthm} in practice.

\section{Preliminaries on Harmonic Maass Forms}

Recall that a \textit{harmonic Maass form} of integer weight $k$ is a real-analytic function $f \colon \mathbb{H} \to \mathbb{C}$ which satisfies the modular transformation law, is annihilated by the weight $k$ hyperbolic Laplacian operator $\Delta_k$, and exhibits at most linear exponential growth at the cusps. If $f$ has poles in $\mathbb{H}$, we say that $f$ is a \textit{polar harmonic Maass form}.

The theory of monstrous moonshine tells us that the Fourier expansion of $j(\tau)$ is the McKay-Thompson series for the identity, meaning that the Fourier coefficients $c(n)$ are the graded dimensions of the Monster module $V^\natural$. From moonshine we also obtain the infinite product identity
$$j(z) - j(\tau) = e^{-2\pi iz}\prod_{m > 0, n \in \mathbb{Z}} \left(1 - e^{2\pi imz}e^{2\pi in\tau}\right)^{c(mn)},$$
known as the \textit{denominator formula} for the Monster Lie algebra.

It turns out that the denominator formula, when viewed as a function of $\tau$, is a polar harmonic Maass form with a simple pole at $z$. More specifically, the denominator formula is equivalent to a theorem of Asai, Kaneko, and Ninomiya (see Theorem 3 of \cite{valuesofmfdivmf}). The theorem states that if we define
$$H_z(\tau) := \sum_{n=0}^{\infty}j_n(z)e^{2\pi in\tau} = \frac{E_4(\tau)^2E_6(\tau)}{\Delta(\tau)} \frac{1}{j(\tau) - j(z)} = -\frac{1}{2\pi i} \frac{j'(\tau)}{j(\tau) - j(z)},$$
then the functions $j_n(\tau)$ form a Hecke system. Namely, if we let $j_0(\tau) = 1$ and $j_1(\tau) = j(\tau) - 744$, then the others are given by
$$j_n(\tau) = j_1(\tau) \mid T(n),$$
where $T(n)$ is the $n$th normalized Hecke operator.

In \cite{divmf} Bringmann et al. generalize the above theorem by constructing weight $2$ polar harmonic Maass forms $H_{N,z}^{\ast}(\tau)$ which generalize the $H_z(\tau)$. Their work extends the result for $j(\tau)$ to meromorphic modular forms on all of the modular curves $X_0(N)$. They also give asymptotics for the Fourier coefficients of the $H_{N,z}^{\ast}(\tau)$ in terms of  ``Ramanujan-like'' expansions, sums of the form
\begin{equation}\label{Ramanujan sums}
\sum_{\substack{\lambda \in \Lambda_z \\ \lambda \leq n}} \sum_{(c,d) \in S_{\lambda}}e\left(-\frac{n}{\lambda}r_z(c,d)\right)e^{\tfrac{2\pi n\Im(z)}{\lambda}}.
\end{equation}
Here we define $e(w) := e^{2\pi iw}$ for $w \in \mathbb{C}$. The definitions of the objects appearing in the sum are as follows. For an arbitrary solution $a,b \in \mathbb{Z}$ to $ad - bc = 1$, we define
\begin{equation*}
\begin{aligned}
r_z(c,d) &:= ac\vert z \vert^2 + (ad + bc)\Re(z) + bd, \\
\Lambda_z &:= \{\alpha\vert z \vert^2 + \beta\Re(z) + \gamma^2 : \alpha, \beta, \gamma \in \mathbb{Z}\}, \\
S_{\lambda} &:= \{(c,d) \in N\mathbb{N}_0 \times \mathbb{Z} : \gcd(c,d) = 1 \textrm{ and } Q_z(c,d) = \lambda\}, \\
Q_z(c,d) &:= c^2\vert z \vert^2 + 2cd\Re(z) + d^2.
\end{aligned}
\end{equation*}
\noindent Note that $r_z(c,d)$ is not uniquely defined. However $e(-nr_z(c,d)/Q_z(c,d))$ is well defined. We quote Theorem 1.1 from \cite{divmf} below.

\begin{thm}\label{PHMF}
If $z \in \mathbb{H}$, then $H^{\ast}_{N,z}(\tau)$ is a weight $2$ polar harmonic Maass form on $\Gamma_0(N)$ which vanishes at all cusps and has a single simple pole at $z$. Moreover, the following are true:
\begin{enumerate}[font=\normalfont]
\item If $z \in \mathbb{H}$ and $\Im(\tau) > \max\{\Im(z), \tfrac{1}{\Im(z)}\}$, then we have that
$$H^{\ast}_{N,z}(\tau) = \frac{3}{\pi[\emph{SL}_2(\mathbb{Z}) : \Gamma_0(N)]\Im(\tau)} + \sum_{n=1}^{\infty}j_{N,n}(z)q^n.$$
\item For $\gcd(N,n) = 1$, we have $j_{N,n}(\tau) = j_{N,1}(\tau) \mid T(n)$.
\item For $n \mid N$, we have $j_{N,n}(\tau) = j_{\tfrac{N}{n},1}(n\tau)$.
\item As $n \to \infty$, we have
$$j_{N,n}(\tau) = \sum_{\substack{\lambda \in \Lambda_{\tau} \\ \lambda \leq n}} \sum_{(c,d) \in S_{\lambda}}e\left(-\frac{n}{\lambda}r_{\tau}(c,d)\right)e^{\tfrac{2\pi n\Im(\tau)}{\lambda}} + O_{\tau}(n).$$
\end{enumerate}
\end{thm}

\noindent If we let $N = 1$, then $j_{1,n}(\tau) = j_n(\tau)$ and we recover the $H_z(\tau)$ up to the addition of the weight $2$ nonholomorphic Eisenstein series $E^{\ast}_2(\tau) := -\tfrac{3}{\pi\Im(\tau)} + E_2(\tau)$.

We also quote Corollary 1.3 from \cite{divmf}, which we will use to obtain the imaginary part of $j^{-1}(\alpha)$.

\begin{cor}\label{Im(z) cor}
Suppose that $f(\tau)$ is a meromorphic modular form of weight $k$ on $\Gamma_0(N)$ whose divisor is not supported at cusps. Let $y_1$ be the largest imaginary part of any points in the divisor of $f(\tau)$ lying in $\mathbb{H}$. Then if $-\frac{1}{2\pi i}\frac{f'(\tau)}{f(\tau)} =: \sum_{n \gg -\infty} a(n)e^{2\pi in\tau}$, we have that
$$y_1 = \lim\sup_{n\to\infty} \frac{\log\vert a(n)\vert}{2\pi n}.$$
\end{cor}

\section{Proof of \Cref{mainthm}}

In this section we will use the results gathered in the previous section to prove \Cref{mainthm}. We first rewrite the asymptotic formula given in \Cref{PHMF} as a sum over corresponding matrices $M = (\begin{smallmatrix} a & b \\ c & d \end{smallmatrix}) \in \Gamma_{\infty} \backslash \Gamma_0(N)$. Direct substitution and simplification shows that $r_z(c,d)/Q_z(c,d) = \Re(Mz)$ and $\Im(z)/Q_z(c,d) = \Im(Mz)$, thus \Cref{PHMF} (4) is equivalent to
$$j_{N,n}(z) = \sum_{\substack{M \in \Gamma_{\infty} \backslash \Gamma_0(N) \\ n\Im(Mz) \geq \Im(z)}}e^{-2\pi inMz} + O_z(n).$$

In the case where $N = 1$, we have $\Im(Mz) \leq \Im(z)$ for all $M \in \textrm{SL}_2(\mathbb{Z})$, thus the $\lambda = 1$ terms dominate in the formula given in \Cref{PHMF} (4). Separating out the $c = 0$ term, we have
$$j_n(z) \approx e^{-2\pi inz} + \sum_{c \geq 1} \sum_{\substack{d \in \mathbb{Z} \\ \gcd(c,d) = 1 \\ \vert cz + d \vert^2 = 1}} e\left(n\frac{d - a}{c}\right)e^{2\pi in\bar{z}},$$
where the $e^{2\pi inz}$ arises from the $c = 0$ term. If $z \in X_0(1) \setminus \{\infty\}$ and $\vert z \vert > 1$, then $Q_z(c,d) = 1$ has no solutions for $c \geq 1$. If $\vert z \vert = 1$ and $z \neq e^{2\pi i/3}$, then the only solution is $(c,d) = (1,0)$ and the second term reduces to $e^{2\pi in\bar{z}}$. Writing $z = x + iy$, we conclude that

\begin{equation}\label{Re(z) asymptotics}
j_n(z)e^{-2\pi ny} \sim \begin{cases}
e^{-2\pi nx} \qquad &\textrm{if $\vert z \vert > 1$} \\
2\cos(2\pi nx) \qquad &\textrm{if $\vert z \vert = 1, z \neq e^{2\pi i/3}$}
\end{cases}.
\end{equation}

\noindent We remark that $j(e^{2\pi i/3}) = 0$ is well known, thus we can exclude the case where $z = e^{2\pi i/3}$.

Let $c(n)$ and $w_n$ be defined as in \Cref{mainthm}. The conditions on $\vert z \vert$ in \Cref{Re(z) asymptotics} are equivalent to the conditions on $\lim_{n\to\infty}\vert a(n) \vert e^{-2\pi ny}$ in the definition of $c(n)$. We will now prove the claimed formulas for $x$ and $y$.

\Cref{Im(z) cor} proves the claimed formula in \Cref{mainthm} for $y$. Once $y \approx y_0$ has been approximated to sufficient precision, we substitute $y_0$ into the formula for $c(n)$. Since the sequence $c(n)$ is bounded, by taking real parts in the case where $\vert z \vert > 1$, it suffices to show the clamed formula for $x$ in the case where $\vert z \vert = 1$.

We have
$$c(n) \approx \cos(w_{n-1} \pm 2\pi x).$$
Now $\vert x \vert \leq \tfrac{1}{2}$ and $w_n \in [0,\pi]$, thus  $w_n \pm 2\pi x \in [-\pi, 2\pi]$. Note that
$$\cos^{-1}(\cos(x_0)) =
\begin{cases}
-x_0 \quad &x_0 \in [-\pi,0) \\
x_0  \quad &x_0 \in [0,\pi] \\
2\pi - x_0 \quad &x_0 \in (\pi,2\pi]
\end{cases}.$$
We thus have
$$\pm w_n \approx w_{n-1} \pm 2\pi x$$
or
$$w_n \approx 2\pi - (w_{n-1} \pm 2\pi x).$$
Rearranging gives the formulas claimed in \Cref{mainthm}.

\section{Examples}

In this section we provide some examples of calculating $z$ using \Cref{mainthm} for selected values of $\alpha$. Throughout this section we let $z = x + iy$ and $q := e^{2\pi in\tau}$. We let $a(n), b(n), c(n), w_n, y_0$ be defined as in \Cref{mainthm}.

\begin{ex} ($\alpha = 2 \cdot 30^3, z = \sqrt{3}i$)

\noindent We have
$$-\frac{1}{2\pi i}\frac{j'(\tau)}{j(\tau) - 2 \cdot 30^3} = 1 + 53256q + 2835807768q^2 + 151013228757024q^3 + \cdots.$$
We find that $b(3) = 1.73205083\ldots$ matches the limiting value up to $7$ decimal places. We see from the size of $y_0$ that $\vert z \vert > 1$. We compute
\begin{equation*}
\begin{aligned}
c(1) &= 1.00007\ldots, \\
c(2) &= 1.00000\ldots.
\end{aligned}
\end{equation*}
We see that $c(n) \to 1$, thus $x = 0$.
\end{ex}

\begin{ex} ($\alpha = -640320^3, z = \tfrac{-1 + \sqrt{163}i}{2}$)

\noindent We have
$$-\frac{1}{2\pi i}\frac{j'(\tau)}{j(\tau) + 640320^3} = 1 - 262537412640768744q + \cdots.$$
We find that $b(1) = 6.3835726674\ldots$ matches the limiting value up to $30$ decimal places. We see from the size of $y_0$ that $\vert z \vert > 1$. We compute
\begin{equation*}
\begin{aligned}
c(1) &= -1.000000000000000000000000000003\ldots, \\
c(2) &= 1.000000000000000000000000000000\ldots.
\end{aligned}
\end{equation*}
We see that $c(n) \approx (-1)^n$, thus $x = -\tfrac{1}{2}$.
\end{ex}

\begin{ex} ($\alpha = 1728, z = i$)

\noindent We have
$$-\frac{1}{2\pi i}\frac{j'(\tau)}{j(\tau) - 1728} = 1 + 984q + 574488q^2 + 307081056q^3 + \cdots.$$
We find that $b(5000) = 1.0000220635600152652\ldots$ matches the limiting value up to $4$ decimal places. We compute
\begin{equation*}
\begin{aligned}
c(1) &= 1.001440\ldots, \\
c(2) &= 0.999503\dots,
\end{aligned}
\end{equation*}
thus $x = 0$.
\end{ex}

\begin{ex} ($\alpha = 1 + i$)

\noindent We have
$$-\frac{1}{2\pi i}\frac{j'(\tau)}{j(\tau) - (1 + i)} = 1 - (744 - i)q + (158280 - 1486i)q^2 - (35797022 - 1065494i)q^3 + \cdots.$$
We find that $b(100) = .8882136152\ldots$. Since $\alpha$ is nonreal we must have $\vert z \vert > 1$. We compute
\begin{equation*}
\begin{aligned}
x &\approx -\frac{1}{2\pi}(w_{103} + w_{102}) = -0.477227209285886\ldots,
\end{aligned}
\end{equation*}
thus
$$z \approx -.4772 + .8882i.$$
To verify the value of $x$, we check that
$$j(-.4772 + .8882i) \approx 1.0042 + 0.9983i.$$
\end{ex}

\section{Acknowledgements}

This research was carried out in fulfillment of the requirements for the M.S. degree in mathematics at Emory University. The author would like to thank Ken Ono for his mentorship and support.

\bibliographystyle{plain}
\bibliography{AnAlgorithmForNumericallyInvertingTheModularJfunctionEalwais}

\begin{thebibliography}{1}

\bibitem{MSBS}
E.~Alwaise.
\newblock An algorithm for numerically computing preimages of the
  $j$-invariant.
\newblock Master's thesis, Emory University, 2017.

\bibitem{valuesofmfdivmf}
J.~Bruinier{,}~W. Kohnen{,} and K.~Ono.
\newblock The arithmetic of the values of modular functions and the divisors of
  modular forms.
\newblock {\em Compositio Mathematica Compositio Math.}, 140(03):552–566,
  2004.

\bibitem{divmf}
K.~Bringmann{,} B. Kane{,} S. L{\"o}brich{,}~K. Ono{,} and L.~Rolen.
\newblock On divisors of modular forms.
\newblock preprint, {arXiv.org:1609.08100v2{,}} 2016.

\end{thebibliography}
\end{document}